\documentclass{amsart}
\usepackage{amssymb, amsmath}
\newcommand\field[1]{\mathbb{#1}}
\newcommand{\CC}{\field{C}}
\newcommand{\NN}{\field{N}}

\newcommand{\Bb}{\mathcal B}

\newcommand{\Hh}{\mathcal H}

\newcommand{\Tt}{\mathcal T}

\newcommand\lsp{\operatorname{span}\nolimits}
\newcommand\clsp{\overline{\lsp}}
\newcommand\proj{\operatorname{proj}}

\newcommand\apbdry[1]{\partial{#1}^{\mathrm ap}}
\newcommand\Pap{P^{\mathrm{ap}}}
\newcommand\Sap{S^{\mathrm{ap}}}

\theoremstyle{plain}
\newtheorem{theorem}{Theorem}[section]

\newtheorem{lemma}[theorem]{Lemma}
\newtheorem{prop}[theorem]{Proposition}
\theoremstyle{definition}

\newtheorem{dfn}[theorem]{Definition}

\numberwithin{equation}{section}

\title[Co-universal algebras of directed graphs]{A direct approach to co-universal algebras associated to directed graphs}

\author[A. Sims]{Aidan Sims}
\address{Aidan Sims \\
School of Mathematics and Applied Statistics \\
Austin Keane building (15) \\
University of Wollongong \\
NSW 2522 \\
AUSTRALIA}
\email{asims@uow.edu.au}

\author[S.B.G. Webster]{Samuel B.G. Webster}
\address{Samuel Webster\\
School of Mathematics and Applied Statistics\\
University of Wollongong\\
NSW 2522\\
AUSTRALIA}
\email{sbgwebster@gmail.com}

\thanks{This research was supported by the Australian Research Council.}
\date{December 10, 2009}

\begin{document}

\maketitle

\begin{abstract}
We prove directly that if $E$ is a directed graph in which
every cycle has an entrance, then there exists a $C^*$-algebra
which is co-universal for Toeplitz-Cuntz-Krieger $E$-families.
In particular, our proof does not invoke ideal-structure theory
for graph algebras, nor does it involve use of the gauge action
or its fixed point algebra.
\end{abstract}

\section{Introduction}

In recent years there has been a great deal of interest in
graph $C^*$-algebras and their generalisations (see
\cite{CBMSbk} for a survey). To associate $C^*$-algebras to a
given generalisation of directed graphs, one assigns partial
isometries to the edges of the graph in a way which encodes
connectivity in the graph. One then aims to identify relations
amongst the partial isometries so that the $C^*$-algebra
universal for these relations satisfies a version of the
Cuntz-Krieger uniqueness theorem. For directed graphs, this
theorem states that if every cycle has an entrance, then any
representation of its Cuntz-Krieger algebra which is nonzero on
generators is faithful. In trying to identify appropriate
relations, there is typically some analogue of a left-regular
representation which points to a natural notion of an abstract
representation; in the case of directed graphs, each graph can
be represented on the Hilbert space with orthonormal basis
indexed by finite paths in the graph, and this gives rise to
the notion of a Toeplitz-Cuntz-Krieger family. However, the
universal $C^*$-algebra for such representations is typically
too big to satisfy a version of the Cuntz-Krieger uniqueness
theorem, and one has to identify an additional relation to
correct this.

In \cite{Katsura2006c}, in the much more general context of
Cuntz-Pimsner algebras associated to Hilbert bimodules, Katsura
developed a very elegant solution to this problem. For directed
graphs, his results say that given any Toeplitz-Cuntz-Krieger
$E$-family consisting of nonzero partial isometries, if the
$C^*$-algebra $B$ it generates carries a circle action
compatible with the canonical gauge action on the Toeplitz
algebra, then there is a canonical homomorphism from $B$ onto
the Cuntz-Krieger algebra. We call this a co-universal
property: the Cuntz-Krieger algebra is co-universal for
Toeplitz-Cuntz-Krieger families consisting of nonzero partial
isometries and compatible with the gauge action. Katsura proved
this theorem \emph{a posteriori}: the Cuntz-Krieger algebra had
already been defined in terms of a universal property. However,
he pointed out that this theorem implies that the Cuntz-Krieger
algebra could be \emph{defined} to be the algebra co-universal
for nonzero Cuntz-Krieger families; one would then have to work
to prove that such a co-universal algebra exists.

When every cycle in the graph $E$ has an entrance, it is a
consequence of the Toeplitz-Cuntz-Krieger uniqueness theorem
that every Cuntz-Krieger $E$-family consisting of nonzero
partial isometries is automatically compatible with the gauge
action. In particular, in this case the use of the gauge action
in Katsura's analysis should not be necessary. In this article
we show that this is indeed the case: we present a direct
argument that for an arbitrary directed graph in which every
cycle has an entrance, there exists a $C^*$-algebra which is
co-universal for Toeplitz-Cuntz-Krieger families consisting of
nonzero partial isometries. In particular, we do not first
identify the Cuntz-Krieger relation or the corresponding
universal $C^*$-algebra. We also do not proceed via the
machinery of ideal structure of Toeplitz algebras of directed
graphs, and we do not deal with the gauge-action of the circle
or with an analysis of its fixed-point algebra. Instead we work
directly with the abelian subalgebra generated by range
projections associated to paths in the graph.

\section{Preliminaries}\label{sec:prelims}

A directed graph $E = (E^0, E^1, r, s)$ consists of a countable
set $E^0$ of vertices, a countable set $E^1$ of edges, and maps
$r,s : E^1 \to E^0$ indicating the direction of the edges. We
will follow the conventions of \cite{CBMSbk} so that a path is
a sequence $\alpha = \alpha_1 \alpha_2 \dots \alpha_n$ of edges
such that $s(\alpha_i) = r(\alpha_{i+1})$ for all $i$. We write
$|\alpha|$ for $n$, and if we want to indicate a segment of a
path, we shall denote it $\alpha_{[p,q]} = \alpha_{p+1}
\alpha_{p+2} \dots \alpha_q$. If $n = \infty$ (so that $\alpha$
is actually a right-infinite string), then we call $\alpha$ an
\emph{infinite path}. The range of an infinite path $\alpha =
\alpha_1\alpha_2\dots$ is $r(\alpha) := r(\alpha_1)$.

For $n \in \NN$, we write $E^n$ for the collection of paths of
length $n$. We write $E^*$ for the category of finite paths in
$E$ (we regard vertices as paths of length zero), and
$E^\infty$ for the collection of all infinite paths. Given
$\alpha \in E^*$ and $X \subset E^* \cup E^\infty$, we denote
$\{\alpha\mu : \mu \in X, r(\mu) = s(\alpha)\}$ by $\alpha X$.
Given $v \in E^0$, a set $X \subset vE^*$ is said to be
\emph{exhaustive} if, for every $\lambda \in vE^*$ there exists
$\alpha \in X$ such that either $\lambda = \alpha\lambda'$ or
$\alpha = \lambda\alpha'$.

Given a directed graph $E$, a \emph{Toeplitz-Cuntz-Krieger
$E$-family} in a $C^*$-algebra $B$ is a pair $(t,q)$ where $t :
e \mapsto t_e$ assigns to each edge a partial isometry in $B$,
and $q : v \mapsto q_v$ assigns to each vertex a projection in
$B$ such that
\begin{itemize}
\item[(TCK1)] the $q_v$ are mutually orthogonal
\item[(TCK2)] each $t_e^* t_e = q_{s(e)}$, and
\item[(TCK3)] for each $v \in E^0$ and each finite subset $F
    \subset vE^1$, we have $q_v \ge \sum_{e \in F} t_e
    t^*_e$.
\end{itemize}
There is a $C^*$-algebra $\Tt C^*(E)$ generated by a
Toeplitz-Cuntz-Krieger family $(s, p)$ which is universal in
the sense that given any Toeplitz-Cuntz-Krieger family $(t,q)$
in a $C^*$-algebra $B$ there is a homomorphism $\pi_{t,q} : \Tt
C^*(E) \to B$ such that $\pi_{t,q}(s_e) = t_e$ and
$\pi_{t,q}(p_v) = q_v$ for all $e \in E^1$ and $v \in E^0$.

Given a path $\alpha \in E^*$ and a Toeplitz-Cuntz-Krieger
$E$-family $(t,q)$, we write $t_\alpha$ for
$t_{\alpha_1}t_{\alpha_2} \dots t_{\alpha_{|\alpha|}}$.

\section{The co-universal algebra}\label{sec:CKUT}

\begin{theorem}\label{thm:main}
Let $E$ be a directed graph, and suppose that every cycle in
$E$ has an entrance. There exists a Toeplitz-Cuntz-Krieger
$E$-family $(\Sap, \Pap)$ consisting of nonzero partial
isometries such that $C^*_{\textrm{min}}(E) := C^*(\Sap, \Pap)$
is co-universal in the sense that given any other
Toeplitz-Cuntz-Krieger $E$ family $(t,q)$ in which each $q_v$
is nonzero, there is a homomorphism $\psi_{t,q} : C^*(t,q) \to
C^*_{\textrm{min}}(E)$ such that $\psi_{t,q}(t_e) =
S^{\mathrm{ap}}_e$ and $\psi_{t,q}(q_v) = \Pap_v$ for all $e
\in E^1$ and $v \in E^0$.
\end{theorem}

Our proof relies on understanding the structure the
$C^*$-algebra generated by the projections $\{t_\lambda
t^*_\lambda : \lambda \in E^*\}$ for a Toeplitz-Cuntz-Krieger
$E$-family $(t,q)$. We begin with the following definition.

\begin{dfn}
Let $E$ be a directed graph. A \emph{boolean representation} of
$E$ in a $C^*$-algebra $B$ is a map $p : \lambda \mapsto
p_\lambda$ from $E^*$ to $B$ such that each $p_\lambda$ is a
projection, and
\[
p_\mu p_\nu
    = \begin{cases}
        p_\nu & \text{ if $\nu = \mu\nu'$} \\
        p_\mu & \text{ if $\mu = \nu\mu'$} \\
        0 & \text{ otherwise.}
    \end{cases}
\]

\end{dfn}

If $p$ is a boolean representation of $E$, then the $p_\lambda$
commute, so $\clsp\{p_\lambda : \lambda \in E^*\}$ is a
commutative $C^*$-algebra.

\begin{lemma}\label{p nu = sum q nu}
Let $E$ be a directed graph, let $p$ be a boolean
representation of $E$, and fix a finite subset $F\subset E^*$.
For $\mu \in F$, define
\[
    q_\mu^F:= p_\mu \prod_{\mu \mu' \in F\setminus\{\mu\}}(p_\mu - p_{\mu\mu'}).
\]
Then the $q_\mu^F$ are mutually orthogonal projections and for
each $\mu \in E^*$,
\begin{equation}\label{eq_p_nu= sum q_nu}
p_\mu = \sum_{\mu\mu' \in F} q_{\mu\mu'}^F
\end{equation}
\begin{proof}
We proceed by induction on $|F|$. If $|F|=1$ then
\eqref{eq_p_nu= sum q_nu} is trivial. Suppose \eqref{eq_p_nu=
sum q_nu} holds whenever $|F|<n$, and fix $F$ with $|F|=n$. Let
$\lambda \in F$ be of maximal length, and let $G = F \setminus
\{\lambda\}$. Then $q_\lambda^F = p_\lambda$, and for $\mu \in
G$,
\[
q_\mu^F
    = \begin{cases}
        q_\mu^G &\text{if } \lambda \neq \mu\mu'\\
        q_\mu^G - q_\mu^G p_\lambda &\text{if } \lambda = \mu\mu'.
    \end{cases}
\]

Fix $\mu \in F$. If $\lambda \neq \mu\mu'$, then the inductive
hypothesis implies that $\sum_{\mu\mu' \in F} q_{\mu\mu'}^F =
\sum_{\mu\mu' \in G} q_{\mu\mu'}^G = p_\mu$. If $\lambda =
\mu\mu'$, then
\begin{flalign*}
    &&\sum_{\mu\mu' \in F} q_{\mu\mu'}^F &= \sum_{\mu\mu' \in G} (q_\mu^G - q_\mu^G p_\lambda) + q_\lambda^F&\\
    &&&= \sum_{\mu\mu' \in G} q_\mu^G - \sum_{\mu\mu' \in G} q_\mu^G p_\lambda + p_\lambda&\\
    &&&= p_\mu -  p_\mu p_\lambda + p_\lambda \qquad\text{by the inductive hypothesis}&\\
    &&&= p_\mu -  p_\lambda + p_\lambda \qquad\text{since $\lambda = \mu\mu'$}&\\
    &&&= p_\mu. &\qedhere
\end{flalign*}
\end{proof}
\end{lemma}

Given a directed graph $E$, we define the set of
\emph{aperiodic boundary paths} in $E$ by
\[
    \apbdry{E} := \big\{\lambda \in E^* : |s(\lambda)E^1| \in \{0,\infty\}\big\}
        \; \sqcup \;E^\infty \setminus \big\{\lambda\mu^\infty : s(\lambda) = r(\mu) = s(\mu)\big\}.
\]
Let $\Hh := \ell^2(\apbdry{E})$, with orthonormal basis
$\{\xi_x : x \in \apbdry{E}\}$, and define $\{\Pap_\lambda :
\lambda \in E^*\} \subset \Bb(\Hh)$ by\label{pg:apboolean}
\[
    \Pap_\lambda \xi_x = \begin{cases} \xi_x &\text{if } x\in\lambda\apbdry{E}\\0 &\text{otherwise.}\end{cases}
\]
Since each $\Pap_\lambda = \operatorname{proj}_{\clsp\{\xi_x :
x \in \lambda\apbdry{E}\}}$, it is straightforward to check
that $\Pap$ is a boolean representation of $E$.

Suppose that every cycle in $E$ has an entrance. We claim that
each $\Pap_\lambda$ is nonzero. Indeed, fix $\lambda \in E^*$.
Since every cycle in $E$ has an entrance, there exists an $x
\in \apbdry{E}$ with $r(x) = s(\lambda)$. Then $\lambda x \in
\apbdry{E}$ so $\Pap_\lambda \xi_{\lambda x} = \xi_{\lambda x}
\neq 0$.

\begin{lemma}\label{faux CK4}
Let $E$ be a directed graph. Let $\lambda \in E^*$, and suppose
that $F \subset s(\lambda) E^*$ is finite and exhaustive. Then
\[
    \prod_{\mu \in F}(\Pap_\lambda - \Pap_{\lambda\mu}) = 0.
\]
\begin{proof}
Let $x \in \apbdry{E}$. We seek $\mu \in F$ such that
$(\Pap_\lambda - \Pap_{\lambda\mu})\xi_x = 0$.

If $\Pap_\lambda \xi_x = 0$, then any $\mu \in F$ will suffice,
so suppose that $\Pap_\lambda \xi_x \neq 0$. Then
$x(0,|\lambda|) = \lambda$. If there exists $\mu \in F$ such
that $x = \lambda\mu x'$, then $(\Pap_\lambda -
\Pap_{\lambda\mu})\xi_x = \xi_x - \xi_x = 0$, so we suppose
that $x \not= \lambda\mu x'$ for all $\mu \in F$ and seek a
contradiction. Since $F$ is exhaustive, there exists $\mu \in
F$ such that $\lambda\mu = x \mu'$ for some $\mu'$ of nonzero
length. In particular, $x$ is a finite path and since $\mu'$
has nonzero length, $s(x)E^1 \not=\emptyset$, so $x \in
\apbdry{E}$ forces $|s(x)E^1| = \infty$. But since no initial
segment of $x$ belongs to $\lambda F$, that $F$ is exhaustive
implies that for each $e \in s(x)E^1$, there exists $\mu_e \in
F$ such that $\lambda\mu_e = xe x_e'$ for some $x_e \in E^*$,
and this contradicts that $F$ is a finite set.
\end{proof}
\end{lemma}

\begin{lemma}\label{converse of CK4}
Let $E$ be a directed graph and $p$ be a boolean representation
of $E$ such that $p_\lambda \neq 0$ for each $\lambda \in E^*$.
If $\{\alpha' : \alpha\alpha' \in F\}$ is not exhaustive, then
$q_\alpha^F \neq 0.$
\begin{proof}
Since $\{\alpha' : \alpha\alpha' \in F\}$ is not exhaustive,
there exists $\tau \in E^*$ such that $\tau
\not=\alpha'\alpha''$ and $\alpha' \not= \tau\tau'$ for each
$\alpha'$ such that $\alpha\alpha' \in F$. In particular, each
$p_{\alpha\alpha'}p_{\alpha\tau} = 0$, and hence $q_\alpha^F
p_{\alpha\tau} = p_{\alpha \tau} \neq 0$, whence $q_\alpha^F
\neq 0.$
\end{proof}
\end{lemma}

\begin{prop}\label{prp:diag cu}
Let $E$ be a directed graph and let $p$ be a boolean
representation of $E$ such that $p_\lambda \neq 0$ for each
$\lambda \in E^*$. Then there is a homomorphism $\psi_p :
\clsp\{p_\lambda : \lambda \in E^*\} \to \clsp\{\Pap_\lambda :
\lambda \in E^*\}$ satisfying $\psi_p(p_\lambda) =
\Pap_\lambda$ for all $\lambda \in E^*$. Moreover, $\psi_p$ is
injective if and only if $\prod_{\mu \in F}(p_\lambda -
p_{\lambda\mu}) = 0$ for all $\lambda \in E^*$ and finite
exhaustive $F \subset s(\lambda)E^*$.
\begin{proof}
For the first assertion it suffices to show that for every
finite subset $F$ of $E^*$ and every collection of scalars
$\{a_\lambda : \lambda \in F\}$,
\begin{equation}\label{Pap to p is norm decr}
    \Big\| \sum_{\lambda \in F} a_\lambda \Pap_\lambda \Big\| \leq \Big\| \sum_{\lambda \in F} a_\lambda p_\lambda\Big\|
\end{equation}
Fix a finite subset $F \subset E^*$, and for $\alpha \in F$,
define $Q^F_\alpha := \Pap_\alpha \prod_{\alpha\alpha' \in F
\setminus \{\alpha\}} (\Pap_\alpha - \Pap_{\alpha\alpha'})$.
For $\lambda \in E^*$, Lemma~\ref{p nu = sum q nu} gives
\[
    \sum_{\lambda \in F} a_\lambda \Pap_\lambda = \sum_{\alpha \in F} \Big(\sum_{\substack{\mu \in F \\ \alpha = \mu\mu'}} a_\mu\Big) Q_\alpha^F,
        \qquad\text{and}\qquad
    \sum_{\lambda \in F} a_\lambda p_\lambda = \sum_{\alpha \in F} \Big(\sum_{\substack{\mu \in F \\ \alpha = \mu\mu'}} a_\mu\Big) q_\alpha^F.
\]
By Lemmas \ref{faux CK4} and \ref{converse of CK4}, we have
$\{\alpha \in F : Q_\alpha^F \neq 0\} \subset \{\alpha \in F :
q_\alpha^F \neq 0\}.$ Hence
\[
\Big\| \sum_{\lambda \in F} a_\lambda \Pap_\lambda \Big\|
    = \max_{Q^F_\alpha \not = 0} \Big| \sum_{\substack{\mu \in F \\ \alpha = \mu\mu'}} a_\mu \Big|
    \le \max_{q^F_\alpha \not = 0} \Big| \sum_{\substack{\mu \in F \\ \alpha = \mu\mu'}} a_\mu \Big|
    = \Big\| \sum_{\lambda \in F} a_\lambda p_\lambda \Big\|,
\]
and the first assertion follows.

For the second assertion, note that if $\prod_{\mu \in
F}(p_\lambda - p_{\lambda\mu}) = 0$ for all $\lambda \in E^*$
and finite exhaustive $F \subset s(\lambda)E^*$, then for each
finite exhaustive $F$, we have $\{\alpha \in F : Q_\alpha^F
\neq 0\} = \{\alpha \in F : q_\alpha^F \neq 0\}$ so the
calculation above shows that $\psi_p$ is isometric.
\end{proof}
\end{prop}

We now show that the $C^*$-algebra generated by any
Toeplitz-Cuntz-Krieger $E$-family in which all the partial
isometries are nonzero admits a conditional expectation onto
the subalgebra spanned by the range projections $s_\lambda
s^*_\lambda$. Recall that given a Toeplitz-Cuntz-Krieger
$E$-family $(t,q)$, we write $\pi_{t,q}$ for the canonical
homomorphism from $\Tt C^*(E)$ to $C^*(t,q)$ induced by the
universal property of the former. We first need a technical
lemma.

\begin{lemma}\label{L rho is a loop}
Suppose $\lambda,\mu,\nu \in E$ satisfy $|\lambda| \geq |\nu| >
|\mu|$ and
\begin{equation}\label{eq:product nonzero}
    t_\lambda t_\lambda^* t_\mu t_\nu^* t_\lambda t_\lambda^* \neq 0.
\end{equation}
Then $\lambda = \nu\nu' = \mu\mu'\nu'$ for some $\mu',\nu'$,
and
\[
    t_\lambda t_\lambda^* t_\mu t_\nu^* t_\lambda t_\lambda^* = t_\lambda t_{\lambda \rho}^*
\]
for some cycle $\rho \in E$.
\end{lemma}
\begin{proof}
Equation~\ref{eq:product nonzero} forces $t_\lambda t^*_\lambda
t_\mu t^*_\mu \not=0$ and $t_\nu t^*_\nu t_\lambda t^*_\lambda
\not= 0$. Hence $\lambda = \nu\nu'$ and $\lambda = \mu \alpha$
for some $\nu'$ and $\alpha$, and since $|\mu| < |\nu|$, this
forces $\nu = \mu\mu'$ and hence $\lambda = \mu\mu'\nu'$ for
some $\mu'$.

We have
\begin{equation}\label{eq rho is a loop}
 0 \neq t_\lambda t_\lambda^* t_\mu t_\nu^* t_\lambda t_\lambda^*
    = t_\lambda t_\lambda^* t_\mu (t_{\mu'}^* t_\mu^*) (t_\mu t_{\mu'} t_{\nu'}) t_\lambda^*
    = t_\lambda t_\lambda^* t_\mu t_{\nu'} t_\lambda^*,
\end{equation}
forcing $s(\mu) = r(\nu')$. Since $r(\mu') = s(\mu)$ and
$s(\mu') = r(\nu')$, $\mu'$ is a cycle, and has nonzero length
since $|\nu| > |\mu|$. Furthermore, continuing from \eqref{eq
rho is a loop}
\[
0 \neq t_\lambda t_\lambda^* t_\mu t_{\nu'} t_\lambda^*
    = (t_\mu t_{\mu'\nu'}) (t_{\mu'\nu'}^* t_\mu^*) t_\mu t_{\nu'} (t_{\nu'}^* t_{\mu\mu'}^*)
        = t_\mu t_{\mu'\nu'} t_{\mu'\nu'}^* t_{\nu'} t_{\nu'}^* t_{\mu\mu'}^*,
\]
so $t^*_{\mu'\nu'} t_\nu'$ is nonzero, forcing
\begin{equation}\label{eq:munu=nurho}
    \mu' \nu' = \nu' \rho\qquad\text{for some $\rho \in E$.}
\end{equation}

We claim $\rho$ is a cycle. We proceed by induction on
$|\nu'|$. As a base case, suppose that $|\nu'| \le |\mu'|$.
Then $\nu' = \mu'_{[0,|\nu'|]}$, so
\[
    \mu' \nu' = \nu' \rho \implies \rho = \mu'_{[|\nu'|,|\mu'|]} \mu'_{[0,|\nu'|]},
\]
whence $r(\rho) = s(\rho) = s(\nu')$.

Now suppose as an inductive hypothesis that $\rho$ is a cycle
whenever~\eqref{eq:munu=nurho} holds with $|\nu'| < n$ for some
$n > |\mu'|$, and fix $\nu'$ with $|\nu'| = n$
satisfying~\eqref{eq:munu=nurho}. In particular, $|\nu'| >
|\mu'|$, so $\mu' =\nu'_{[0,|\mu'|]}$ and $\mu' \nu'_{[|\mu'|,
|\nu'|]} = \nu' = \nu'_{[|\mu'|, |\nu'|]} \rho$. Since
$|\nu'_{[|\mu'|, |\nu'|]}| = |\nu'| - |\mu'| < n$, the
inductive hypothesis now implies that $\rho$ is a cycle.

Finally, from \eqref{eq rho is a loop},
\begin{flalign*}
  &&0 \not= t_\lambda t_\lambda^* t_\mu t_\nu^* t_\lambda t_\lambda^*
    &= t_\lambda t_\lambda^* t_\mu t_{\nu'} t_\lambda^* &\\
    &&&= t_\lambda (t_{\mu'\nu'}^* t_\mu^*) t_\mu t_{\nu'} t_\lambda^*
     = t_\lambda t_{\nu' \rho}^* t_{\nu'} t_\lambda^*
     = t_\lambda t_\rho^* t_{\nu'}^* t_{\nu'} t_\lambda^* = t_\lambda t_\rho^* t_\lambda^*,&\qedhere
\end{flalign*}
\end{proof}

Given a directed graph $E$ and $e \in E^1$, define $\Sap_e \in \Bb(\ell^2(\apbdry{E}))$ by
\[
\Sap_e \xi_x
    = \begin{cases}
        \xi_{ex} &\text{ if $s(e) = r(x)$}\\
        0 &\text{ otherwise.}
    \end{cases}
\]
Then with $\{\Pap_v : v \in E^0\}$ as on
page~\pageref{pg:apboolean}, $(\Sap,\Pap)$ is a
Toeplitz-Cuntz-Krieger $E$-family, and $\Sap_\lambda
(\Sap_\lambda)^* = \Pap_\lambda$ for all $\lambda \in E^*$. We
denote $C^*(\Sap,\Pap)$ by $C^*_{\mathrm{min}}(E)$.

\begin{prop}\label{prp:expectation}
Let $E$ be a directed graph, and suppose that every cycle in
$E$ has an entrance. Let $(t,q)$ be a Toeplitz-Cuntz-Krieger
$E$-family. Then $q : \lambda \mapsto t_\lambda t^*_\lambda$ is
a boolean representation of $E$, and there exists a conditional
expectation $\Phi_{t,q} : C^*(t,q) \to \clsp\{q_\lambda :
\lambda \in E^*\}$ satisfying $\Phi_{t,q}(t_\mu t^*_\nu) =
\delta_{\mu,\nu} q_\mu$. In particular we have
\begin{equation}\label{eq:compatible}
\psi_{q} \circ \Phi_{t,q} \circ \pi_{t,q} = \Phi_{\Sap,\Pap} \circ \pi_{\Sap,\Pap}.
\end{equation}
\end{prop}
\begin{proof}
It is standard that the $q_\lambda$ form a boolean
representation of $E$.

\textbf{Claim~1.} Given a finite subset $F$ of $E^*$ and a
collection $\{a_{\mu,\nu} : \mu,\nu \in F\}$ of scalars,
\[
\Big\|\sum_{\mu \in F} a_{\mu,\mu} q_\mu\Big\| \le \Big\|\sum_{\mu,\nu \in F} a_{\mu,\nu} t_\mu t^*_\nu\Big\|.
\]
By enlarging $F$ (and setting the extra scalars equal to zero),
we may assume that $F$ is closed under initial segments in the
sense that if $\mu\nu \in F$ then $\mu \in F$.

For each $\lambda \in F$, let $T^F_\lambda := \{\lambda' \in
s(\lambda)E^* : \lambda\lambda' \in F, |\lambda'| > 0\}$. For
each $\lambda \in F$ such that $T^F_\lambda$ is not exhaustive,
fix a path $\alpha^\lambda$ such that $\alpha^\lambda \not=
\mu\mu'$ and $\mu \not= \alpha^\lambda\alpha'$ for all $\mu \in
T^F_\lambda$. Since every cycle in $E$ has an entrance,
\cite[Lemma~3.7]{CBMSbk} implies that for each $v$ such that $v
= s(\alpha^\lambda)$ for some $\lambda$, there exists $\tau^v
\in v E^*$ such that either: (1) $s(\tau^v) E^1 = \emptyset$;
or (2) $|\tau^v|
> \max\{|\lambda| : \lambda \in F, T^F_\lambda\text{ is not exhaustive}\}$,
and $\tau^v_k \not= \tau^v_{|\tau^v|}$ for all $k < |\tau^v|$.
We write $\tau^\lambda$ for $\tau^{s(\alpha^\lambda)}$.

For each $\lambda \in F$ we define
\[
\phi^F_\lambda := \begin{cases}
        q_{\lambda\alpha^\lambda\tau^\lambda} &\text{ if $T^F_\lambda$ is not exhaustive} \\
        q^F_\lambda &\text{ otherwise}.
    \end{cases}
\]
By definition we have each $\phi^F_\lambda \le q^F_\lambda$,
and since the $q^F_\lambda$ are mutually orthogonal, it follows
that the $\phi^F_\lambda$ are also. Hence
\[
\Big\|\sum_{\mu,\nu \in F} a_{\mu,\nu} t_\mu t^*_\nu\Big\|
    \ge \Big\|\sum_{\lambda \in F} \phi^F_\lambda \Big(\sum_{\mu,\nu \in F} a_{\mu,\nu} t_\mu t^*_\nu\Big) \phi^F_\lambda \Big\|.
\]

Fix $\lambda,\mu,\nu \in F$. We claim that
\begin{equation}\label{eq:comp by phi}
\phi^F_\lambda t_\mu t^*_\nu \phi^F_\lambda
    = \begin{cases}
        \phi^F_\lambda &\text{ if $\mu=\nu$ and $\lambda = \mu\lambda'$ for some $\lambda'$} \\
        0 &\text{ otherwise.}
    \end{cases}
\end{equation}
To see this, suppose first that $\mu = \nu$. If $\lambda =
\mu\lambda'$ then $\phi^F_\lambda \le t_\lambda t^*_\lambda \le
t_\mu t^*_\mu$ by definition of $\phi^F_\lambda$. If $\lambda
\not= \mu\lambda'$, then either $\mu = \lambda\mu'$ in which
case $\phi^F_\lambda \le (t_\lambda t^*_\lambda -
t_{\lambda\mu'} t^*_{\lambda\mu'}) \perp t_\mu t^*_\mu$, or
else $\mu \not= \lambda\mu'$ in which case $\phi^F_\lambda \le
t_\lambda t^*_\lambda \perp t_\mu t^*_\mu$.

Now suppose that $\mu \not= \nu$; by symmetry under adjoints,
we may assume that $|\mu| < |\nu|$. We must show that
$\phi^F_\lambda t_\mu t^*_\nu \phi^F_\lambda = 0$. Since
$\phi^F_\lambda \le t_\lambda t^*_\lambda$, if $t_\lambda
t_\lambda^* t_\mu t_\nu^* t_\lambda t_\lambda^* = 0$ then we
are done, so we may assume that $t_\lambda t_\lambda^* t_\mu
t_\nu^* t_\lambda t_\lambda^* \not= 0$. Then Lemma~\ref{L rho
is a loop} implies that $\lambda = \nu\nu' = \mu\mu'\nu'$ and
that $t_\lambda t_\lambda^* t_\mu t_\nu^* t_\lambda t_\lambda^*
= t_\lambda t_{\lambda \rho}^*$ for some cycle $\rho \in E$.
Hence $\phi^F_\lambda \le t_\lambda t^*_\lambda$ forces
\[
\phi_\lambda^F t_\mu t^*_\nu \phi_\lambda^F
    = \phi_\lambda^F t_\lambda t^*_{\lambda\rho} \phi_\lambda^F.
\]
We consider two cases: $T^F_\lambda$ is exhaustive, or it is
not.

First suppose that $T^F_\lambda$ is exhaustive. Fix $n$ such
that $n|\rho| > \max\{|\lambda'| : \lambda' \in T^F_\lambda\}$.
Then $\rho^n = \lambda' \beta$ for some $\lambda' \in
T_\lambda^F$ and $\beta \in E^*$. In particular, $\lambda' =
\rho^n_{[0,|\lambda'|]},$ forcing $\rho_1 = \lambda'_1$. Since
$\lambda\lambda' \in F$ and $F$ is closed under initial
segments, $\lambda \rho_1 \in F$ and then
\[
t_{\lambda\rho}t^*_{\lambda\rho} \phi_\lambda^F
    \le t^*_{\lambda \rho} t^*_{\lambda \rho} (t_\lambda t^*_\lambda - t_{\lambda \rho_1}t^*_{\lambda \rho_1}) = 0,
\]
giving $\phi^F_\lambda s_\mu s^*_\nu \phi^F_\lambda = 0$.

Now suppose that $T^F_\lambda$ is not exhaustive. Then
$\phi^F_\lambda = q_{\lambda\alpha^\lambda\tau^\lambda}$. We
then have
\[
\phi_\lambda^F t_\mu t^*_\nu \phi_\lambda^F
    = t_{\lambda\alpha^\lambda\tau^\lambda} t^*_{\mu'\nu'\alpha^\lambda\tau^\lambda}
        t_{\nu'\alpha^\lambda\tau^\lambda} t^*_{\lambda\alpha^\lambda\tau^\lambda}.
\]
This is nonzero only if $\mu'\nu'\alpha^\lambda\tau^\lambda =
\nu'\alpha^\lambda\tau^\lambda\zeta$ for some $\zeta$. This is
impossible if $s(\tau^\lambda)E^1 = \emptyset$, so suppose that
$s(\tau^\lambda)E^1 \not= \emptyset$. By choice of
$\tau^\lambda$, we have $|\tau^\lambda|
> |\mu'|$. Let $m = |\nu'\alpha^\lambda\tau^\lambda|$. Then
\[
(\mu'\nu'\alpha^\lambda\tau^\lambda)_m
    = \tau^\lambda_{|\tau^\lambda| - |\mu'|}
    \not= \tau^\lambda_{|\tau^\lambda|}
    = (\nu'\alpha^\lambda\tau^\lambda)_m,
\]
so $\mu'\nu'\alpha^\lambda\tau^\lambda \not=
\nu'\alpha^\lambda\tau^\lambda\zeta$ for all $\zeta$, and hence
$\phi^F_\lambda t_\mu t^*_\nu \phi^F_\lambda = 0$,
establishing~\eqref{eq:comp by phi}.

We now have
\begin{align*}
\Big\|\sum_{\mu \in F} a_{\mu,\mu} q_\mu\Big\|
    &= \Big\|\sum_{\mu \in F} \Big(\sum_{n \le |\mu|} a_{\mu_{[0,n]},\mu_{[0,n]}}\Big) q^F_\mu \Big\| \\
    &= \max_{\mu \in F} \Big|\sum_{n \le |\mu|} a_{\mu_{[0,n]},\mu_{[0,n]}}\Big| \\
    &= \Big\|\sum_{\mu \in F} \Big(\sum_{n \le |\mu|} a_{\mu_{[0,n]},\mu_{[0,n]}}\Big) \phi^F_\mu \Big\| \\
    &= \Big\|\sum_{\lambda \in F} \phi^F_\lambda \Big(\sum_{\mu,\nu\in F} a_{\mu,\nu} t_\mu t^*_\nu\Big) \phi^F_\lambda\Big\|
        \qquad\text{by~\eqref{eq:comp by phi}}\\
    &\le \Big\|\sum_{\mu,\nu \in F} a_{\mu,\nu} t_\mu t^*_\nu\Big\|
\end{align*}
completing the proof of Claim~1.

Claim~1 implies that the formula $t_\mu t^*_\nu \mapsto
\delta_{\mu,\nu} t_\mu t^*_\mu$ extends to a well-defined
linear map $\Phi_{t,q}$ from $C^*(t,q)$ to $\clsp\{q_\lambda :
\lambda \in E^*\}$. This $\Phi_{t,q}$ is a linear idempotent of
norm 1, and hence a conditional expectation (see for example
\cite[Definition~II.6.10.2 and
Theorem~II.6.10.2]{Blackadar2006}).

The final statement is straightforward since the two maps in
question agree on spanning elements.
\end{proof}


\begin{lemma}\label{lem:cu exp faithful}
Let $E$ be a directed graph in which every cycle has an
entrance. Then the expectation $\Phi_{\Sap,\Pap} :
C^*_{\mathrm{min}}(E) \to \clsp\{\Pap_\lambda : \lambda \in
E^*\}$ obtained from Proposition~\ref{prp:expectation} is
faithful on positive elements.
\end{lemma}
\begin{proof}
It suffices to show that for $a \in C^*_{\mathrm{min}}(E)$,
$\Phi_{\Sap,\Pap}(a)$ is equal to the strong-operator limit
$\sum_{x \in \apbdry{E}} \big(a \xi_x | \xi_x) \theta_{\xi_x,
\xi_x}$ for all $a \in C^*_{\mathrm{min}}(E)$, where the
$\xi_x$ are the canonical orthonormal basis for
$\ell^2(\apbdry{E})$ and $\theta_{\xi_x, \xi_x}$ is the
rank-one projection onto $\CC \xi_x$. Fix $\mu,\nu \in E^*$ and
$x \in \apbdry{E}$. We have
\[
\big(\Sap_\mu (\Sap_\nu)^* \xi_x \big| \xi_x\big)
    = \big((\Sap_\nu)^* \xi_x \big| \big (\Sap_\mu)^*\xi_x\big)
    = \begin{cases}
        1 &\text{ if $x = \nu y = \mu y$ for some $y$}\\
        0 &\text{ otherwise.}
    \end{cases}
\]
Since $x \in \apbdry{E}$, we have $y \not= \rho^\infty$ for any
cycle $\rho$, so $\mu y = \nu y$ forces $|\mu| = |\nu|$, and
hence $\mu = \nu$.

Hence
\[
\sum_{x \in \apbdry{E}}
\big(\Sap_\mu (\Sap_\nu)^* \xi_x | \xi_x) \theta_{\xi_x, \xi_x}
    = \delta_{\mu,\nu} \proj_{\clsp\{\xi_x : x \in \mu\apbdry{E}\}}
    = \delta_{\mu,\nu} \Pap_\mu
\]
as required.
\end{proof}

We now have the tools we need to prove the main theorem.

\begin{proof}[Proof of Theorem~\ref{thm:main}]
We showed that the $\Pap_\lambda$, and in particular the
$\Pap_v$ are nonzero immediately subsequent to their
definition.

Fix a Toeplitz-Cuntz-Krieger $E$-family $(t,q)$ with each $q_v$
nonzero. We will show that $\ker(\pi_{t,q}) \subset
\ker(\pi_{\Sap, \Pap})$ in $\Tt C^*(E)$, and hence that
$\pi_{\Sap, \Pap}$ descends to the desired homomorphism
$\psi_{t,q} : C^*(t,q) \to C^*_{\mathrm{min}}(E)$.

Since each $t_\lambda^* t_\lambda = q_{s(\lambda)}$, each
$q_\lambda$ is nonzero, so Proposition~\ref{prp:diag cu}
implies that there is a homomorphism $\psi_q : \clsp\{q_\lambda
: \lambda \in E^*\} \to \clsp\{\Pap_\lambda : \lambda \in
E^*\}$ taking each $q_\lambda$ to $\Pap_\lambda$.

We calculate
\begin{flalign}
&&\pi_{t,q}(a) = 0
    &\iff \pi_{t,q}(a^*a) = 0 \nonumber & \\
    &&&\implies \psi_q \circ \Phi_{t,q} \circ \pi_{t,q}(a^*a) = 0 \label{eq:key implic} & \\
    &&&\iff \Phi_{\Sap,\Pap}(\pi_{\Sap,\Pap}(a^*a)) = 0 \quad\text{ by~\eqref{eq:compatible}} \nonumber & \\
    &&&\iff \pi_{\Sap,\Pap}(a^*a) \quad\text{ by Lemma~\ref{lem:cu exp faithful}} \nonumber & \\
    &&&\iff \pi_{\Sap,\Pap}(a) = 0. \nonumber &\qedhere
\end{flalign}
\end{proof}

It is, of course, interesting to know when the homomorphism
$\psi_{t,q}$ of Theorem~\ref{thm:main} is injective.

\begin{theorem}
Let $E$ be a directed graph in which every cycle has an
entrance.
\begin{enumerate}
\item\label{it:couniv} If $(t,q)$ is a
    Toeplitz-Cuntz-Krieger family with each $q_v$ nonzero,
    then the homomorphism $\psi_{t,q}$ of
    Theorem~\ref{thm:main} is injective if and only if (a)
    $\prod_{\lambda \in F} (q_v - q_\lambda) = 0$ whenever
    $v \in E^0$ and $F \subset vE^*$ is finite exhaustive;
    and (b) the expectation $\Phi_{t,q}$ is faithful.
\item\label{it:injective} If $\pi$ is homomorphism from
    $C^*_{\mathrm{min}}(E)$ to a $C^*$-algebra $C$ such
    that each $\pi_{\Pap_v}$ is nonzero, then $\pi$ is
    injective.
\end{enumerate}
\end{theorem}
\begin{proof}
(\ref{it:couniv}) Since conditions (a)~and~(b) hold in
$C^*_{\mathrm{min}}(E)$, the ``only if" implication is trivial.
For the ``if" implication, note that given a
Toeplitz-Cuntz-Krieger $E$-family $(t,q)$, we have
$\ker(\pi_{t,q}) = \ker(\pi_{\Sap, \Pap})$ whenever the
implication~\eqref{eq:key implic} is equivalence. Condition~(a)
implies that $\psi_q$ is faithful by the final statement of
Proposition~\ref{prp:diag cu}, and this combined with~(b)
implies that $\psi_q \circ \Phi_{t,q}$ is faithful on positive
elements, giving
\[
\pi_{t,q}(a^*a) = 0 \iff \psi_q \circ \Phi_{t,q} \circ \pi_{t,q}(a^*a) = 0
\]
as required.

(\ref{it:injective}) Define a Toeplitz-Cuntz-Krieger $E$-family
by $t_e := \pi(\Sap_e)$ and $q_v := \pi(\Pap_v)$.
Theorem~\ref{thm:main} supplies a homomorphism $\psi_{t,q} :
C^*(t,q) \to C^*_{\mathrm{min}}(E)$ which is an inverse for
$\pi$.
\end{proof}

\end{document}